\RequirePackage{fix-cm}
\documentclass{svjour3}                     % onecolumn (standard format)
\smartqed  % flush right qed marks, e.g. at end of proof
\usepackage{mathbbold}
\usepackage{mathrsfs}
\usepackage{makeidx,ifthen}
\usepackage{latexsym}
\usepackage{graphicx}
\usepackage{amsmath}
\usepackage{amssymb}
\usepackage{amsfonts}
\usepackage{wasysym}

\usepackage{mathptmx}      % use Times fonts if available on your TeX system
\newtheorem{thm}{Theorem}%%[section]

\newtheorem{lem}[thm]{Lemma}
\newtheorem{prop}[thm]{Proposition}

\numberwithin{equation}{section}

\begin{document}

\title{Solutions for a class of quasilinear
Schr\"{o}dinger equations with critical exponents term
\thanks{Supported by NSFC under grant numbers 11201488 and 11101418.
}
}


\author{Zhouxin Li         \and
        Yimin Zhang %etc.
}


\institute{Z. Li (Corresponding Author)\at
              Department of Mathematics, Central South University, Changsha, China {\rm 400083} \\
              %Tel.: +86-731-89920238\\
              \email{lzx@math.pku.edu.cn}           %  \\
           \and
           Y. Zhang \at
              Wuhan Institute of Physics and Mathematics, Chinese Academy of Sciences, \\Wuhan, China {\rm 430071}\\
              \email{zhangym802@126.com}
}

 \date{}%%Received: date / Accepted: date}

\maketitle

\begin{abstract}
In this paper,
we study a class of quasilinear Schr\"{o}dinger equation of the form
\begin{eqnarray*}
-\varepsilon^2\Delta u+V(x)u-\varepsilon^2(\Delta(|u|^{2\alpha}))|u|^{2\alpha-2}u
    &=&\lambda|u|^{q-2}u+|u|^{2^*(2\alpha)-2}u,\quad\mbox{in }{\mathbb{R}}^N,
\end{eqnarray*}
where $\varepsilon>0,~\lambda>0,~q\geq2,~\alpha>1/2$ are constants,
$N\geq3$. By using change of variable and variational approach,
 the existence of positive solution which has a local maximum point and decays exponentially is obtained.
$\\${\bf Mathematics Subject Classification (2000)} {35J10 \and 35J20 \and 35J25}
\end{abstract}

\section{Introduction}
\label{sec:1}
Let $l$ and $h$ be real functions of pure power forms,
it is interesting to consider the existence of solutions to
the following quasilinear Schr\"{o}dinger equation
\begin{equation}\label{sch1}
  i\partial_tz=-\varepsilon^2\Delta z+W(x)z-l(|z|^2)z-k\varepsilon^2\Delta h(|z|^2)h'(|z|^2)z,\quad x\in{\mathbb{R}}^N,N\geq3,
\end{equation}
where $W(x)$ is a given potential, $k$ is a real constant,
$\varepsilon>0$ is a real parameter. It has many applications in
physics according to different types of $h$. For example, in
\cite{Ku-81}, it was used in plasma physics with $h(s)=s$, and was
used in \cite{Ri-94} to models the self-channeling of high-power
ultrashort laser in matter with $h(s)=(1+s)^{1/2}$. Readers can
refer to \cite{Liu-Wang-03,Liu-Wang2-03} for references of more
applications of it.

In this paper, we assume that $h(s)=s^{\alpha}$, $l(s)=\lambda
s^{(q-2)/2}+s^{2^*(\alpha)-1}$, where
$\lambda>0,~q\geq2,~\alpha>1/2,~2^*=2N/(N-2)$ are constants. If we
consider standing waves solutions of the form
$z(x,t)=\exp(-iEt/\varepsilon)u(x)$, then $z(x,t)$ satisfies
equation (\ref{sch1}) if and only if the function $u(x)$ solves the
equation
\begin{eqnarray}\label{sch2}
-\varepsilon^2\Delta u&+&V(x)u-k\alpha\varepsilon^2(\Delta(|u|^{2\alpha}))|u|^{2\alpha-2}u \nonumber\\
    &=&\lambda|u|^{q-2}u+|u|^{2^*(2\alpha)-2}u,\quad u>0,x\in{\mathbb{R}}^N,
\end{eqnarray}
where $V(x)=W(x)-E$ is the new potential function.

In case $k=0$, equation (\ref{sch2}) is a semilinear elliptic
equation which has been extensively studied in the past two decades.
In recent years, the quasilinear case $k\neq0$ arose a lot of
interest to mathematical researchers. From \cite{Liu-Wang2-03}, we
know that the constant $2^*(2\alpha)>2^*$ is thought to behave like
the critical exponent to equation (\ref{sch2}). When $\alpha=1$
(i.e. $h(s)=s$) and $\varepsilon=1$, this kind of problems with
different types of nonlinearities $l(s)$ at sub-"critical growth",
i.e. at sub-$2^*(2\alpha)$ growth, have been widely studied, see
\cite{Ch-Gu-05,Co-Je-04,Do-Se-09,Liu-Wang2-03} and so on. In
\cite{Liu-Wang2-03}, by using a changing of variable, they transform
the equation to a semilinear one, then the existence of solutions
was obtained via variational methods under different types of
potentials $V(x)$. This method is significant and was widely used in
the studies of this kind of problems. For general $\alpha$, there
are few results for this case, as far as we know, just \cite{aw,aw2,Liu-Wang-03,Mo-06-2}.
In \cite{Liu-Wang-03}, for
$\alpha>\frac{1}{2}$, $l(s)=\lambda s^{\frac{p-1}{2}}$ and
$2<p+1<2^*(2\alpha)$, the existence of solution of equation
(\ref{sch2}) without critical term was obtained by using the method
of Lagrange multiplier. In \cite{Mo-06-2}, the existence of at least
one or sometimes two standing wave solutions for
$\alpha>\frac{1}{2}$ and $l(s)=\mu f(x) s^{\frac{p-1}{2}}$ was
obtained through fibreing method. Employing the change variable method
just as \cite{Liu-Wang2-03}, the authors of \cite{aw} obtained
the existence of at least one positive solution for
$\alpha>\frac{1}{2}$ and general $l(s)$ by using variational
approached. Moreover, in \cite{aw2}, for $V(x)=\lambda$,
$\alpha>\frac{1}{2}$ and $l(s)=s^{\frac{p-1}{2}}$, they obtained the
unique existence of positive radial solution under some suitable
conditions.

Problem at "critical growth", i.e. at $2^*(2\alpha)$ growth rate,
was studied by Moameni \cite{Mo-06} with a general $\varepsilon>0$
and $\alpha=1$. It was assumed in \cite{Mo-06} that $V(x)=0$ on an
annule, which enable the avoidance of proving the compact embedding
near the origin. In \cite{llw}, the existence of radially symmetric solution was
obtained for $\varepsilon>0$ small enough. Other kinds of such
problems at "critical growth" were studied in \cite{do-Mi-So-07,llw,Zh-Wa-Sh} and the references
   therein. But all these are studied for $\alpha=1$, for general
   $\alpha$, there is no results according to what we know.

In this paper, our aim is to study the existence of positive
solutions of (\ref{sch2}) with general $\alpha>1/2$ and at
$2^*(2\alpha)$ growth. Problem of (\ref{sch2}) at $2^*(2\alpha)$
growth has two difficulties. Firstly, the embedding
$H^1({\mathbb{R}}^N)\hookrightarrow
L^{2^*(2\alpha)}({\mathbb{R}}^N)$ is not compact, so it is hard to
prove the Palais-Smale ((PS) in short) condition. Secondly, even if
we can obtain the compactness result of (PS) sequence, it is only
holds at some level of positive upper bound, it is difficult for us
to prove that the functional has such minimax level.

We assume that $V(x)$ is locally H\"{o}lder continuous and
\begin{description}
  \item[(V)] $\exists~ V_\infty>V_0>0$ such that $\min_{x\in{\mathbb{R}}^N}V(x)=V_0$ and $\lim_{|x|\to\infty}V(x)=V_\infty$.
\end{description}

Note that assumption (V) allows zero be the minimum point of $V(x)$.
This is different to the assumptions on $V(x)$ in \cite{Mo-06}.

Under assumption (V), we define a space
\[
X:=\bigg\{u\in H^{1}({\mathbb{R}}^N):\int_{{\mathbb{R}}^N}V(x)u^2<\infty\bigg\}
\]
with the norm
$
\|u\|_X^2=\int_{{\mathbb{R}}^N}|\nabla u|^2+\int_{{\mathbb{R}}^N}V(x)u^2.
$

For simplicity of notation, we let $m=2\alpha$, $\bar q=q/m$ and
$k\alpha=1$ in this paper. Set
\[
g(t)=\lambda|t|^{q-2}t+|t|^{2^*m-2}t\quad\mbox{and}\quad G(t)=\int_0^tg(s)ds.
\]
We formulate problem (\ref{sch2}) in variational structure in the space $X$ as follows:
\[
I(u)=\frac{\varepsilon^2}{2}\int_{{\mathbb{R}}^N}(1+m|u|^{2(m-1)})|\nabla u|^2+\frac{1}{2}\int_{{\mathbb{R}}^N}V(x)u^2-\int_{{\mathbb{R}}^N}G(u)
\]
Note that $I$ is lower simicontinuous on $X$, we define that $u\in X$ is a {\it weak solution} for (\ref{sch2}) if
$u\in X\cap L^{\infty}({\mathbb{R}}^N)$ and it is a critical point of $I$.

Firstly, for an arbitrary $\varepsilon>0$, we have:

\begin{thm} \label{main1}
Assume that $q\in(2m,2^*m)$ and that condition (V) holds.

$\\$
Case 1: $1< m<2$. Assume that one of the following conditions holds:
\begin{description}
  \item[(i)] $\bar q>4+\frac{2}{m}$ if $N=3$;
  \item[(ii)] $\bar q>2+\frac{2}{m}$ if $N=4$;
  \item[(iii)] $\bar q>\frac{4}{3}+\frac{2}{m}$ if $N=5$;
  \item[(iv)] $\bar q>2$ if $N\geq6$.
\end{description}
Case 2: $m\geq2$. Assume that one of the following conditions holds:
\begin{description}
  \item[(v)] $\bar q>5$ if $N=3$;
  \item[(vi)] $\bar q>3$ if $N=4$;
  \item[(vii)] $\bar q>\frac{7}{3}$ if $N=5$;
  \item[(viii)] $\bar q>2$ if $N\geq6$.
\end{description}
Then for $\varepsilon>0$ small enough,
problem (\ref{sch2}) has a positive weak solution $u_\varepsilon\in X\cap L^{\infty}({\mathbb{R}}^N)$ with
\[
\lim_{\varepsilon\to0}\|u_\varepsilon\|_{X}=0,\quad\mbox{and}\quad
u_\varepsilon(x)\leq C\exp(-\frac{\beta}{\varepsilon}|x-x_\varepsilon|).
\]
where  $C>0,\beta>0$ are constants,
$x_\varepsilon\in{\mathbb{R}}^N$ is a local maximum point of $u_\varepsilon$.
\end{thm}

Next, we consider the case $\varepsilon=1$.
We have the following result:

\begin{thm} \label{main2}
Assume that all conditions in Theorem \ref{main1} hold and that $\varepsilon=1$, then
problem (\ref{sch2}) has a positive weak solution $u_1\in X\cap L^{\infty}({\mathbb{R}}^N)$.
\end{thm}

This paper is organized as follows. In section 2, we first use a change of variable
to reformulate the problem, then we modify the functional in order to
regain the (PS) condition. In section 3, we prove that the functional
satisfies the (PS) condition, this is a crucial job of this paper. Finally, in section 4,
we prove the main theorems, which involves the construction of a mountain pass level at a certain high.

\section{ Preliminaries}
\label{sec:2}

Since $I$ is lower semicontinuous on $X$, we follow the idea in
\cite{Co-Je-04,Liu-Wang2-03} and make the change of variables
$v=f^{-1}(u)$, where $f$ is defined by
\begin{eqnarray*}
\left\{
  \begin{array}{ll}
    f(0)=0, & \hbox{} \\
    f'(v)=(1+m|f(v)|^{2(m-1)})^{-1/2}, & \hbox{on $[0,+\infty)$;} \\
    f(v)=-f(-v), & \hbox{on $(-\infty,0]$.}
  \end{array}
\right.
\end{eqnarray*}
The above function $f(t)$ and its derivative satisfy the following
properties (see \cite{aw,aw2,Liu-Wang2-03}):

\begin{lem}\label{lem.f}
For $m>1$, we have

(1) $f$ is uniquely defined, $C^2$ and invertible;

(2) $|f'(t)|\leq1$ for all $t\in{\mathbb{R}}$;

(3) $|f(t)|\leq |t|$ for all $t\in{\mathbb{R}}$;

(4) $f(t)/t\to1$ as $t\to0$;

(5) $|f(t)|\leq m^{1/2m}|t|^{1/m}$  for all $t\in{\mathbb{R}}$;

(6) $\frac{1}{m}f(t)\leq tf'(t)\leq f(t)$  for all $t>0$;

(7) $f(t)/\sqrt[m]{t}\to m^{1/2m}$ as $t\to+\infty$;
\end{lem}

According to \cite{Do-Se-09} (see Corollary 2.1 and Proposition 2.2 in it,
note that the embedding in Corollary 2.1 of \cite{Do-Se-09} is also compact.), we have:
\begin{lem}\label{lem.imbed}
The map: $v\mapsto f(v)$ from $X$ into $L^r({\mathbb{R}}^N)$ is continuous for
$1\leq r\leq 2^*m$, and is locally compact for $1\leq r< 2^*m$.
\end{lem}

Using this change of variable, we rewrite the functional $I(u)$ to:
\[
J(v)=I(f(v))=\frac{\varepsilon^2}{2}\int_{{\mathbb{R}}^N}|\nabla v|^2
+\frac{1}{2}\int_{{\mathbb{R}}^N}V(x)f^2(v)-\int_{{\mathbb{R}}^N}G(f(v)).
\]
The critical point of $J$ is the weak solution of equation
\begin{eqnarray}\label{Sch3}
-\varepsilon^2\Delta v+V(x)f(v)f'(v)=g(f(v))f'(v),\quad  x\in{\mathbb{R}}^N.
\end{eqnarray}

Now we define a suitable modification of the functional $J$ in order
to regain the Palais-Smale condition. In this time, we make use of
the method in \cite{Pi-Fe-96}.

Let $l$ be a positive constant such that
\begin{equation}\label{def.l}
l=\sup\{s>0: \frac{g(t)}{t}\leq\frac{V_0}{k}\mbox{ for every }0\leq t\leq s\}
\end{equation}
for some $k>\theta/(\theta-2)$ with $\theta\in(2m,q]$. We define the functions:
\[
\gamma(s)=\left\{
            \begin{array}{ll}
              g(s), & \hbox{$s>0$;} \\
              0, & \hbox{$s\leq0$.}
            \end{array}
          \right.\quad
\bar\gamma(s)=\left\{
                \begin{array}{ll}
                  \gamma(s), & \hbox{$0\leq s\leq l$;} \\
                  \frac{V_0}{k}s, & \hbox{$s>l$.}
                \end{array}
              \right.
\]
and
\[
p(x,s)=\chi_R(x)\gamma(s)+(1-\chi_R(x))\bar\gamma(s),
\]
\[
P(x,s)=\int_0^sp(x,t)dt,
\]
where $\chi_R$ denotes the characteristic function of the set $B_R$
(the ball centered at 0 and with radius $R$ in ${\mathbb{R}}^N$),
$R>0$ is sufficiently large and such that
$$\min_{B_R}V(x)<\min_{\partial B_R}V(x).$$
By definition, the function $p(x,s)$ is measurable in $x$, of class $C$ in $s$ and satisfies:
\begin{description}
  \item[(p$_1$)] $0<\theta P(x,s)\leq p(x,s)s$ for every $x\in B_R$ and $s\in{\mathbb{R}}^+$.
  \item[(p$_2$)] $0\leq2P(x,s)\leq p(x,s)s\leq\frac{1}{k}V(x)s^2$ for every $x\in B_R^c:={\mathbb{R}}^N\setminus B_R$ and $s\in{\mathbb{R}}^+$.
\end{description}

Now we study the existence of solutions for the deformed equation:
\begin{eqnarray}\label{Sch4}
-\varepsilon^2\Delta v+V(x)f(v)f'(v)=p(x,f(v))f'(v),\quad x\in{\mathbb{R}}^N.
\end{eqnarray}
The corresponding functional of (\ref{Sch4}) is given by
\[
\bar J(v)=\frac{\varepsilon^2}{2}\int_{{\mathbb{R}}^N}|\nabla v|^2
+\frac{1}{2}\int_{{\mathbb{R}}^N}V(x)f^2(v)-\int_{{\mathbb{R}}^N}P(x,f(v)).
\]
For $v\in X$, since
\begin{eqnarray}\label{lem.X.1}
|\nabla (|f(v)|^m)|^2=\frac{m^2|f(v)|^{2(m-1)}}{1+m|f(v)|^{2(m-1)}}|\nabla v|^2\leq m|\nabla v|^2,
\end{eqnarray}
we infer that $|f(v)|^m\in X$. By Sobolev inequality, we have
\begin{eqnarray}\label{lem.X.2}
  \|f(v)\|_{2^*m} = \||f(v)|^{m}\|_{2^*}^{1/m}\leq C\|\nabla(|f(v)|^{m})\|_2^{1/m}\leq C\|v\|_{X}^{1/m}.
\end{eqnarray}
It results that $f(v)\in L^{2^*m}({{\mathbb{R}}^N})$.
Using interpolation inequality, we obtain that $f(v)\in L^{q}({{\mathbb{R}}^N})$.
Thus $\bar J$ is well defined on $X$.
Let $(v_n)\subset X, v\in X$ with $v_n\to v$ in $X$. Then from Lemma \ref{lem.imbed}, we infer that $V(x)f^2(v_n)\to V(x)f^2(v)$ in $L^1({{\mathbb{R}}^N})$ and that $f(v_n)\to f(v)$ in $L^{q}({{\mathbb{R}}^N})$.
Thus $\bar J$ is continuous on $X$.
$\bar J$ is Gateaux-differentiable in $X$ and the G-derivative is
\begin{eqnarray*}
\langle \bar J'(v),\varphi\rangle&=&\varepsilon^2\int_{{\mathbb{R}}^N}\nabla v\nabla \varphi+
        \int_{{\mathbb{R}}^N}V(x)f(v)f'(v)\varphi \\
    &&-\int_{{\mathbb{R}}^N}p(x,f(v))f'(v)\varphi,
        \quad\forall\varphi\in X.
\end{eqnarray*}
Then if $v\in X\cap L^\infty({\mathbb{R}}^N)$ is a critical point of $\bar J$,
and $v(x)\leq a:=f^{-1}(l),~\forall x\in B_R^c$,
we have $u=f(v)\in X\cap L^\infty({\mathbb{R}}^N)$ (note that we have $|u|\leq|v|$ and $|\nabla u|\leq|\nabla v|$ by the properties of $f$) is a solution of (\ref{sch2}).

\section{ Compactness of (PS) sequence}
\label{sec:3}

In this section, we show that the functional $\bar J$
satisfies (PS) condition, this is a crucial job, its proof is
composed of four steps. Let $S$ denotes the best Sobolev constant,
we have

\begin{lem}\label{lem.PS}
Assume that condition (V) holds and $q\in(2m,2^*m)$. Then $\bar J$
satisfies (PS) condition at level $c_\varepsilon<\frac{1}{Nm}\varepsilon^N S^{N/2}$.
\end{lem}
{\bf Proof}\quad Let $(v_n)\in E$ be a (PS) sequence of $\bar J$ at level $c_\varepsilon$, that is, $(v_n)$ satisfies:
\begin{eqnarray}\label{ps-j-1}
    \bar J(v_n)&=&\frac{\varepsilon^2}{2}\int_{{\mathbb{R}}^N}|\nabla v_n|^2+\frac{1}{2}\int_{{\mathbb{R}}^N}V(x)f^2(v_n) \nonumber\\
    &&-\int_{{\mathbb{R}}^N}P(x,f(v_n))=c_\varepsilon+o(1),
\end{eqnarray}
and
\begin{eqnarray}\label{ps-j-2}
    \langle\bar J'(v_n),\varphi\rangle&=&\varepsilon^2\int_{{\mathbb{R}}^N}\nabla v_n\nabla \varphi+
        \int_{{\mathbb{R}}^N}V(x)f(v_n)f'(v_n)\varphi \nonumber\\
    &&-\int_{{\mathbb{R}}^N}p(x,f(v_n))f'(v_n)\varphi
    =o(1)\|\varphi\|_X,\quad \forall \varphi\in X.
\end{eqnarray}

We divide the proof into four steps.

Step 1: the sequence
$
\int_{{\mathbb{R}}^N}(|\nabla v_n|^2+V(x)f^2(v_n))
$
is bounded.
Multiplying (\ref{ps-j-1}) by $\theta$ ($\theta$ is given in section 2) and using ($p_1$)-($p_2$), we get
\begin{eqnarray*}
  &&\frac{\theta\varepsilon^2}{2}\int_{{\mathbb{R}}^N}|\nabla v_n|^2+\frac{\theta}{2}\int_{{\mathbb{R}}^N}V(x)f^2(v_n)\\
  &&\quad\leq \int_{B_R}p(x,f(v_n))f(v_n)+\frac{\theta}{2k}\int_{B_R^c}V(x)f^2(v_n)+\theta c_\varepsilon+o(1).
\end{eqnarray*}
On the other hand, taking $\varphi=f(v_n)/f'(v_n)$ in (\ref{ps-j-2}), we get
\begin{eqnarray*}
  && \int_{{\mathbb{R}}^N}\varepsilon^2\bigg(1+\frac{m(m-1)|f(v_n)|^{2(m-1)}}{1+m|f(v_n)|^{2(m-1)}}\bigg)|\nabla v_n|^2+\int_{{\mathbb{R}}^N}V(x)f^2(v_n) \\
  &&\quad = \int_{{\mathbb{R}}^N}p(x,f(v_n))f(v_n)+o(\|v_n\|_X)
  \geq \int_{B_R}p(x,f(v_n))f(v_n)+o(1)\|v_n\|_X.
\end{eqnarray*}
Conbining the above two inequalities, we get
\begin{eqnarray}\label{lem2-6-0}
  &&(\frac{\theta}{2}-m)\varepsilon^2\int_{{\mathbb{R}}^N}|\nabla v_n|^2+(\frac{\theta}{2}-\frac{\theta}{2k}-1)\int_{{\mathbb{R}}^N}V(x)f^2(v_n)  \nonumber\\
  &&\quad\leq \frac{\theta\varepsilon^2}{2} \int_{{\mathbb{R}}^N}|\nabla v_n|^2-
    \int_{{\mathbb{R}}^N}\varepsilon^2\bigg(1+\frac{m(m-1)|f(v_n)|^{2(m-1)}}{1+m|f(v_n)|^{2(m-1)}}\bigg)|\nabla v_n|^2 \nonumber\\
  &&\qquad +\frac{\theta}{2}\int_{{\mathbb{R}}^N}V(x)f^2(v_n)-\frac{\theta}{2k}\int_{B_R^c}V(x)f^2(v_n)-\int_{{\mathbb{R}}^N}V(x)f^2(v_n) \nonumber\\
  &&\quad\leq \theta c_\varepsilon+o(1)+o(1)\|v_n\|_X.
\end{eqnarray}
Since $\theta>2m$ and $k>\frac{\theta}{\theta-2}$, we get the conclusion from (\ref{lem2-6-0}).

Step 2: for every $\delta>0$, there exists $R_1\geq R>0$ such that
\begin{eqnarray}\label{lem2-6-1}
\limsup_{n\to\infty}\int_{B_{2R_1}^c}(|\nabla v_n|^2+V(x)f^2(v_n))<\delta.
\end{eqnarray}
We consider a cut-off function $\psi_{R_1}=0$ on $B_{R_1}$,
$\psi_{R_1}=1$ on $B_{2R_1}^c$, $|\nabla \psi_{R_1}|\leq C/R_1$ on ${\mathbb{R}}^N$ for some
constant $C>0$. On one hand, taking $\varphi=f(v_n)/f'(v_n)$,
we compute $\langle\bar J'(v_n),\varphi\psi_{R_1}\rangle$ and get
\begin{eqnarray}\label{lem2-6-2}
  o(1)\|v_n\|_X &=& \int_{{\mathbb{R}}^N}\varepsilon^2\bigg(1+\frac{m(m-1)|f(v_n)|^{2(m-1)}}{1+m|f(v_n)|^{2(m-1)}}\bigg)|\nabla v_n|^2\psi_{R_1} \nonumber\\
  &&+\int_{{\mathbb{R}}^N}\varepsilon^2\varphi\nabla v_n\nabla \psi_{R_1}
  +\int_{{\mathbb{R}}^N}V(x)f^2(v_n)\psi_{R_1} \nonumber\\
  &&-\int_{{\mathbb{R}}^N}p(x,f(v_n))f(v_n)\psi_{R_1} \nonumber\\
  &\geq& \int_{{\mathbb{R}}^N}\varepsilon^2|\nabla v_n|^2\psi_{R_1}
  +\int_{{\mathbb{R}}^N}\varepsilon^2\varphi\nabla v_n\nabla \psi_{R_1} \nonumber\\
  &&+(1-\frac{1}{k})\int_{{\mathbb{R}}^N}V(x)f^2(v_n)\psi_{R_1}.
\end{eqnarray}
On the other hand, by H\"{o}lder inequality,
\begin{eqnarray}\label{lem2-6-3}
\bigg|\int_{{\mathbb{R}}^N}\varphi\nabla v_n\nabla \psi_{R_1}\bigg|\leq\frac{C}{R_1}\|\nabla v_n\|_{L^2({\mathbb{R}}^N)}\|\varphi\|_{L^2({\mathbb{R}}^N)}.
\end{eqnarray}
Note that $\|\nabla v_n\|_{L^2({\mathbb{R}}^N)}$ is bounded, and
\begin{eqnarray}\label{lem2-6-4}
    \|\varphi\|_{L^2({\mathbb{R}}^N)}^2
    &=&\int_{{\mathbb{R}}^N}f^2(v_n)(1+m|f(v_n)|^{2(m-1)}) \nonumber\\
    &=&\int_{{\mathbb{R}}^N}f^2(v_n)+m\int_{{\mathbb{R}}^N}|f(v_n)|^{2m},
\end{eqnarray}
by (\ref{lem.X.2}), $\|\varphi\|_{L^2({\mathbb{R}}^N)}$ is bounded also. Therefore, it follows from (\ref{lem2-6-2})-(\ref{lem2-6-4})
that
\[
\limsup_{n\to\infty}\int_{B_{2R_1}^c}(|\nabla v_n|^2+V(x)f^2(v_n))\leq\frac{C}{R_1}
\]
for $R_1$ sufficiently large, this yields (\ref{lem2-6-1}).

Step 3, there exists $v\in X$ such that
\begin{eqnarray}\label{lem2-6-5}
\lim_{n\to\infty}\int_{{\mathbb{R}}^N}p(x,f(v_n))f(v_n)=\int_{{\mathbb{R}}^N}p(x,f(v))f(v).
\end{eqnarray}
Firstly, by step 1, there exists $v\in X$ such that up to a subsequence, $v_n\to v$ weakly in $X$ and $v_n\to v$ a.e. in ${\mathbb{R}}^N$.
Since we may replace $v_n$ by $|v_n|$, we assume $v_n\geq0$ and $v\geq0$.
By (\ref{lem2-6-1}), for any $\delta>0$, there exists $R_1>0$ sufficiently large such that
\[
\limsup_{n\to\infty}\int_{B_{2R_1}^c}(|\nabla v_n|^2+V(x)f^2(v_n))\leq k\delta.
\]
Therefore, by ($p_2$) we have
\begin{equation}\label{lem2-6-6}
    \limsup_{n\to\infty}\int_{B_{2R_1}^c}p(x,f(v_n))f(v_n)
    \leq \limsup_{n\to\infty}\int_{B_{2R_1}^c}\frac{V(x)}{k}f^2(v_n)\leq\delta,
\end{equation}
and by Fatou Lemma,
\begin{equation}\label{lem2-6-7}
    \int_{B_{2R_1}^c}p(x,f(v))f(v)\leq\delta.
\end{equation}
Secondly, we prove that
\begin{equation}\label{lem2-6-8}
    \int_{B_{2R_1}}p(x,f(v_n))f(v_n)\to\int_{B_{2R_1}}p(x,f(v))f(v).
\end{equation}
Then from this, (\ref{lem2-6-6})-(\ref{lem2-6-7}), and the arbitrariness of $\delta$, we get (\ref{lem2-6-5}).
In fact, since $(v_n)$ is bounded in $X$, we have $(f(v_n))$ is bounded also.
Thus there exists a $w\in X$ such that
    $f(v_n)\rightharpoonup w$ in $X$,
    $f(v_n)\to w$ in $L^r(B_{R_1})$ for $1\leq r<2^*$ and
    $f(v_n)\to w$ a.e. in $B_{R_1}$.
According to (\ref{lem.X.1}), $(|f(v_n)|^m)$ is also bounded in $X$.
By a normal argument, we have
    $|f(v_n)|^m\rightharpoonup |w|^m$ in $X$,
    $|f(v_n)|^m\to |w|^m$ in $L^r(B_{R_1})$ for $1\leq r<2^*$ and
    $|f(v_n)|^m\to |w|^m$ a.e. in $B_{R_1}$.
Applying Lions' concentration compactness principle \cite{PLL-85} to $(|f(v_n)|^m)$ on $\bar B_{R_1}$, we obtain that
there exist two nonnegative measures $\mu,~\nu$,
a countable index set $K$, positive
constants $\{\mu_k\},~ \{\nu_k\},~ k\in K$ and a collection of points
$\{x_k\},~ k\in K$ in $\bar{B}_{R_1}$ such that for all $k\in K$,
\begin{description}
  \item[(i)] $\nu=|w|^{2^*m}+\sum\limits_{k\in K}\nu_k\delta_{x_k}$;
  \item[(ii)] $\mu=|\nabla (|w|^m)|^{2}+\sum\limits_{k\in K}\mu_k\delta_{x_k}$;
  \item[(iii)] $\mu_k\geq S\nu_k^{2/2^*}$,
\end{description}
where $\delta_{x_k}$ is the Dirac measure at $x_k$,
$S$ is the best Sobolev constant.
We claim that $\nu_k=0$ for all $k\in K$.
In fact, let $x_k$ be a singular point of measures $\mu$ and $\nu$,
as in \cite{Li-Zou-98}, we define a function $\phi\in
C_0^{\infty}({\mathbb{R}}^N)$ by
\[
\phi(x)=\left\{
  \begin{array}{ll}
    1, & \hbox{$B_{\rho}(x_k)$;} \\
    0, & \hbox{${\mathbb{R}}^N\backslash B_{2\rho}(x_k)$;} \\
    \phi\geq0,|\nabla\phi|\leq\frac{1}{\rho}, & \hbox{$B_{2\rho}(x_k)\backslash B_{\rho}(x_k)$.}
  \end{array}
\right.
\]
where $B_{\rho}(x_k)$ is a ball centered at $x_k$ and with radius $\rho>0$.
We take $\varphi=\phi f(v_n)/f'(v_n)$ as test functions in $\langle \bar J'(v_n),\varphi\rangle$ and get
\begin{eqnarray}\label{prop4-3-1}
&&\int_{{\mathbb{R}}^N} \varepsilon^2\bigg(1+\frac{m(m-1)|f(v_n)|^{2(m-1)}}{1+m|f(v_n)|^{2(m-1)}}\bigg)|\nabla v_n|^2\cdot\phi\nonumber\\
    &&\quad +\int_{{\mathbb{R}}^N}\varepsilon^2\nabla v_n\nabla\phi\cdot f(v_n)/f'(v_n)
        +\int_{{\mathbb{R}}^N}V(x)f^2(v_n)\phi\nonumber\\
    &&\qquad -\int_{{\mathbb{R}}^N}p(x,f(v_n))f(v_n)\phi=o(1)\|v_n\phi\|_X.
\end{eqnarray}
Then Lions' concentration compactness principle implies that
\begin{eqnarray}\label{prop4-3-2}
\int_{B_{R_1}}|\nabla |f(v_n)|^m|^2\phi\to\int_{B_{R_1}}\phi
d\mu,\quad\int_{B_{R_1}}|f(v_n)|^{2^*m}\phi\to\int_{B_{R_1}}\phi d\nu.
\end{eqnarray}
Since $x_k$ is singular point of $\nu$, by the continuity of $f$,
we have $$f(v_n(x))|_{(B_{2\rho}\setminus\{x_k\})}\to\infty$$ as $\rho\to0$.
Thus
\[
1+\frac{m(m-1)|f(v_n)|^{2(m-1)}}{1+m|f(v_n)|^{2(m-1)}}=m-o(\rho).
\]
on $B_{2\rho}$ for $\rho$ sufficiently small.
Then by (\ref{lem.X.1}) we get from (\ref{prop4-3-1}) that
\begin{eqnarray}\label{prop4-3-3}
&&\int_{B_{R_1}}\varepsilon^2\phi d\mu-\int_{B_{R_1}}\phi d\nu\nonumber\\
&&\quad=\lim_{n\to\infty}\bigg[\int_{B_{R_1}}\varepsilon^2|\nabla |f(v_n)|^m|^2\phi-\int_{B_{R_1}}|f(v_n)|^{2^*m}\phi\bigg]\nonumber\\
&&\quad\leq\lim_{n\to\infty}\bigg[\int_{B_{R_1}} m\varepsilon^2|\nabla v_n|^2\phi-\int_{B_{R_1}}|f(v_n)|^{2^*m}\phi\bigg]\nonumber\\
&&\quad\leq\lim_{n\to\infty}\bigg[\int_{B_{R_1}}
     \varepsilon^2\bigg(1+\frac{m(m-1)|f(v_n)|^{2(m-1)}}{1+m|f(v_n)|^{2(m-1)}}\bigg)|\nabla v_n|^2\phi\nonumber\\
        &&\qquad+o(\rho)\int_{B_{R_1}}\varepsilon^2|\nabla v_n|^2\phi-\int_{B_{R_1}}|f(v_n)|^{2^*m}\phi\bigg]\nonumber\\
&&\quad\leq\lim_{n\to\infty}\bigg[-\int_{B_{R_1}} \varepsilon^2\nabla v_n\nabla \phi\cdot f(v_n)/f'(v_n)
        +\lambda\int_{B_{R_1}} |f(v_n)|^q\phi\nonumber\\
        &&\qquad+o(\rho)\int_{B_{R_1}}\varepsilon^2|\nabla v_n|^2\phi+o(1)\|v_n\phi\|_X\bigg].
\end{eqnarray}
We prove that the last
inequality in (\ref{prop4-3-3}) tends to zero as $\rho\to0$.
By H\"{o}lder inequality, we have
\begin{eqnarray}\label{prop4-3-4}
&&\lim_{n\to\infty}\bigg|\int_{B_{R_1}} \nabla v_n\nabla\phi\cdot f(v_n)/f'(v_n)\bigg|\nonumber\\
&&\quad\leq\limsup_{n\to\infty}\bigg(\int_{B_{R_1}}|\nabla v_n|^2\bigg)^{1/2}\cdot
    \bigg(\int_{B_{R_1}} |[f(v_n)/f'(v_n)]\cdot\nabla\phi|^2\bigg)^{1/2}.
\end{eqnarray}
Since $|f(v_n)/f'(v_n)|^2=f^2(v_n)+m|f(v_n)|^{2m}$, using H\"{o}lder inequality we have
\begin{eqnarray*}
&&\lim_{n\to\infty}\int_{B_{R_1}} |[f(v_n)/f'(v_n)]\cdot\nabla\phi|^2 \\
&&\quad \leq C\rho\big(\|w\|^2_{L^{2^*}(B_{2\rho}(x_j))}+\|w\|^2_{L^{2^*m}(B_{2\rho}(x_j))}\big)
\to0
\end{eqnarray*}
as $\rho\to0$. Thus we obtain that the right hand side of (\ref{prop4-3-4}) tends to 0.
On the other hand, since $q\in(2m,2^*m)$, by Lemma \ref{lem.imbed}, we can prove that
$g(x,h(v_n))h(v_n)\phi\to g(x,w)w\phi$ in $L^1({B_{R_1}})$
and
$\int_{B_{R_1}} g(x,w)w\phi\to0$ as $\rho\to0$.
All these facts imply that the last
inequality in (\ref{prop4-3-3}) tends to zero as $\rho\to0$. Thus $\nu_k\geq\varepsilon^2\mu_k$. This means that either
$\nu_k=0$ or $\nu_k\geq \varepsilon^N S^{N/2}$ by virtue of Lions' concentration compactness principle. We claim that the latter
is impossible. Indeed, if $\nu_k\geq \varepsilon^N S^{N/2}$ holds for
some $k\in K$, then
\begin{eqnarray}
c_\varepsilon&=&\lim_{n\to\infty}\bigg{\{}\bar J(v_n)-\frac{1}{2m}\langle \bar J'(v_n),f(v_n)/f'(v_n)\rangle\bigg{\}}\nonumber\\
&\geq&\lim_{n\to\infty}\bigg{\{}(\frac{1}{2m}-\frac{1}{2^*m})\int_{{\mathbb{R}}^N}|f(v_n)|^{2^*m}\bigg{\}}
\geq(\frac{1}{2m}-\frac{1}{2^*m})\int_{{\mathbb{R}}^N} d\nu \nonumber\\
&\geq&(\frac{1}{2m}-\frac{1}{2^*m})\int_{{\mathbb{R}}^N}|w|^{2^*m}+(\frac{1}{2m}-\frac{1}{2^*m})S^{N/2}\varepsilon^N
\geq\frac{1}{Nm}\varepsilon^N S^{N/2},\nonumber
\end{eqnarray}
which is a contradiction. Thus $\nu_k=0$ for all $k\in K$, and it
implies that $\|f(v_n)\|_{L^{2^*m}(B_{R_1})}$ $\to\|w\|_{L^{2^*m}(B_{R_1})}$. By the uniform convexity
of $L^{2^*m}(B_{R_1})$, we have $f(v_n)\to w$ strongly in
$L^{2^*m}(B_{R_1})$. Finally, since $p(x,f(v_n))f(v_n)$ is sub-$(2^*m)$ growth on $B_{2{R_1}}\setminus B_{R_1}$,
we conclude that (\ref{lem2-6-8}) holds. This proves (\ref{lem2-6-5}).

Step 4: $(v_n)$ is compact in $X$. Since we have (\ref{lem2-6-5}), the proof of the compactness is trivial.
This completes the proof of the lemma.
$\quad\square$

\section{ Proof of main results}
\label{sec:4}

Before we prove Theorem \ref{main1}, we will show firstly some
properties about the change variable $f$.
\begin{lem}\label{lem.f1}
Let $f_1(v)=|f(v)|^m/v,~v\neq0$ and $f_1(0)=0$,
then $f_1$ is continuous, odd, nondecreasing and
\begin{eqnarray}\label{lem.f.limit}
\lim_{v\to0}f_1(v)=0, \quad\mbox{and } \lim_{|v|\to+\infty}|f_1(v)|=\sqrt{m}.
\end{eqnarray}
\end{lem}
{\bf Proof}\quad In fact, by (6) of Lemma \ref{lem.f},
\[
f'_1(v)=v^{-2}(m|f(v)|^{m-2}f(v)f'(v)v-|f(v)|^m)\geq0,
\]
so $f_1$ is nondecreasing. By (4) of Lemma \ref{lem.f}, $f_1(v)\to0$ as $v\to0$. Finally, according to Hospital Principle,
\[
\lim_{v\to+\infty}f_1(v)=\lim_{v\to+\infty}\frac{|f(v)|^m}{v}=\lim_{v\to+\infty}m|f(v)|^{m-2}f(v)f'(v)=\sqrt{m}.
\]
This shows that (\ref{lem.f.limit}) holds.
$\quad\Box$

\begin{lem}\label{lem.f2}
There exists $d_0>0$ such that $$\lim_{v\to+\infty}(\sqrt{m}v-f^m(v))\geq d_0.$$
\end{lem}
{\bf Proof}\quad
Assume that $v>0$. Since by (6) of Lemma \ref{lem.f}, $f(v)\leq mf'(v)v$, we have
\begin{eqnarray}\label{lem.f2.1}
\sqrt{m}v-f^m(v)&\geq&\sqrt{m}v-mf^{m-1}(v)f'(v)v\nonumber\\
    &=&\frac{\sqrt{1+mf^{2(m-1)}(v)}-\sqrt{m}f^{m-1}(v)}{\sqrt{1+mf^{2(m-1)}(v)}}\sqrt{m}v \nonumber\\
    &=&\frac{\sqrt{m}v}{\big(\sqrt{1+mf^{2(m-1)}(v)}+\sqrt{m}f^{m-1}(v)\big)\sqrt{1+mf^{2(m-1)}(v)}} \nonumber\\
    &\geq& \frac{\sqrt{m}v}{2(1+mf^{2(m-1)}(v))}\geq\frac{f^m(v)}{4mf^{2(m-1)}(v)} \nonumber\\
    &=&\frac{1}{4mf^{m-2}(v)}:=d(m,v).
\end{eqnarray}
In the last inequality, we have used the fact that $\sqrt{m}v\geq f^m(v)$
and that $mf^{2(m-1)}(v)>1$ for $v>0$ sufficiently large.

If $1< m< 2$, then $d(m,v)\to+\infty$ as $v\to+\infty$.
If $m=2$, then $d(m,v)=1/8$.
If $m>2$, we claim that $\sqrt{m}v- f^m(v)\to0$ as $v\to+\infty$ is impossible. In fact,
assume on the contrary, then using Hospital Principle, we get
\begin{eqnarray*}
0&\leq& \lim_{v\to+\infty}\frac{\sqrt{m}v- f^m(v)}{f^{2-m}(v)}
    =\lim_{v\to+\infty}\frac{\sqrt{m}-mf^{m-1}(v)f'(v)}{(2-m)f^{1-m}(v)f'(v)}\\
    &=& \lim_{v\to+\infty}\frac{m}{(2-m)f^{1-m}(v)\big(\sqrt{m}\sqrt{1+mf^{2(m-1)}(v)}+mf^{m-1}(v)\big)}\\
    &=&\frac{1}{2(2-m)}<0.
\end{eqnarray*}
This is a contradiction. Thus for all $m>1$, there exists $d_0>0$
such that there holds
$$\lim_{v\to+\infty}(\sqrt{m}v-f^m(v))\geq d_0.$$
This completes the proof.
$\quad\Box$

\begin{lem}\label{lem.f3}
We have

(i) If $1< m <2$, then
\[
\lim_{v\to+\infty}\frac{\sqrt{m}v-f^m(v)}{f^{2-m}(v)}=\frac{1}{2(2-m)}.
\]

(ii) If $m\geq2$, then
\[
\lim_{v\to+\infty}\frac{\sqrt{m}v-f^m(v)}{\log {f(v)}}\leq\left\{
                                                            \begin{array}{ll}
                                                              \frac{1}{2}, & \hbox{$m=2$;} \\
                                                              0, & \hbox{$m>2$.}
                                                            \end{array}
                                                          \right.
\]
\end{lem}
{\bf Proof}\quad Firstly, we prove part (i). According to (\ref{lem.f2.1}) in Lemma \ref{lem.f2}, we have
$\sqrt{m}v-f^m(v)\to+\infty$ as $v\to+\infty$. Thus by Hospital Principle, we get
\[
\lim_{v\to+\infty}\frac{\sqrt{m}v-f^m(v)}{f^{2-m}(v)}
    =\lim_{v\to+\infty}\frac{\sqrt{m}-mf^{m-1}(v)f'(v)}{(2-m)f^{1-m}(v)f'(v)}
    =\frac{1}{2(2-m)}.
\]

Next, we prove part (ii).  If there exists a constant $C>0$ such that $\sqrt{m}v-f^m(v)\leq C$,
then the conclusion holds. Otherwise, assume that $\sqrt{m}v-f^m(v)\to+\infty$ as $v\to+\infty$.
Then again by Hospital Principle, we have
\begin{eqnarray*}
\lim_{v\to+\infty}\frac{\sqrt{m}v-f^m(v)}{\log{f(v)}}
    =\lim_{v\to+\infty}\frac{\sqrt{m}-mf^{m-1}(v)f'(v)}{f'(v)/f(v)}
    =\left\{
          \begin{array}{ll}
            \frac{1}{2}, & \hbox{$m=2$;} \\
            0, & \hbox{$m>2$.}
          \end{array}
        \right.
\end{eqnarray*}
This completes the proof.
$\quad\Box$

To prove Theorem \ref{main1}, it is crucial to prove that $\bar J$ has the mountain pass level
$c_\varepsilon<\frac{1}{Nm}\varepsilon^N S^{N/2}$.
Let us consider the following family of functions in \cite{Br-Ni-83}
\[
v^*_\omega(x)=\frac{[n(n-2)\omega^2]^{(n-2)/4}}{[\omega^2+|x|^2]^{(n-2)/2}},
\]
which solves the equation $-\Delta u=u^{2^*-1}$ in ${\mathbb{R}}^N$ and satisfies
$\|\nabla v^*_\omega\|^2_{L^2}=\|v^*_\omega\|^{2^*}_{L^{2^*}}=S^{N/2}$.
Let $\omega$ be such that $2\omega<R$ and
let $\eta_\omega(x)\in[0,1]$ be a positive smooth cut-off function with
$\eta_\omega(x)=1$ in $B_\omega$, $\eta_\omega(x)=0$ in $B_R\setminus B_{2\omega}$.
Let $v_\omega=\eta_\omega v^*_\omega$. For all $\omega>0$, there exists
$t^\omega>0$ such that $\bar J(t^\omega v_\omega)<0$ for all $t>t^\omega$.
Define the class of paths
\[
\Gamma=\{\gamma\in C([0,1],X):~\gamma(0)=0,\gamma(1)=t^\omega v_\omega\},
\]
and the minimax level
\[
c_\varepsilon=\inf_{\gamma\in\Gamma}\max_{t\in[0,1]}\bar J(\gamma(t)).
\]
Let $t_\omega$ be such that
\[
\bar J(t_\omega v_\omega)=\max_{t\geq0}\bar J(tv_\omega).
\]
Note that the sequence $(v_\omega)$ is uniformly bounded in $X$,
then if $\bar J(t_\omega v_\omega)\to0$ as $t_\omega\to0$, we are done;
on the other hand, if $t_\omega\to+\infty$, then $\bar J(t_\omega v_\omega)\to-\infty$,
which is impossible, so it remains to consider the case where the sequence $(t_\omega)$
is upper and lower bounded by two positive constants.
According to \cite{Br-Ni-83}, we have, as $\omega\to0$,
\[
\|\nabla v_\omega\|^2_{L^{2}}=S^{N/2}+O(\omega^{N-2}),\quad
\|v_\omega\|^{2^*}_{L^{2^*}}=S^{N/2}+O(\omega^{N}).
\]
Let $a\in(0,\frac{\varepsilon^{(N-2)/2}}{2\sqrt{m}})$, $b\in(\frac{2\varepsilon^{(N-2)/2}}{\sqrt{m}},+\infty)$
be such that
$t_\omega\in[a,b],~\forall \omega\in(0,\omega_0)$, where $\omega_0>0$ small enough.
By computing $\frac{d}{dt}\bar J(tv_\omega)=0$,
we obtain $t_\omega=\frac{\varepsilon^{(N-2)/2}}{\sqrt{m}}+o(1)$.
Let
\[
H(v)=-\frac{1}{2}V(x)f^2(v)+\frac{\lambda}{q}|f(v)|^q-\frac{1}{2^*m}|\sqrt{m}v|^{2^*}+\frac{1}{2^*m}|f(v)|^{2^*m},
\]
then by (\ref{lem.f.limit}) and (4) of Lemma \ref{lem.f}, for $m>1$, we have
\[
\lim_{|v|\to+\infty}H(v)/|v|^{2^*}=0,\quad\mbox{and}\quad \lim_{v\to0}H(v)/v^2=-\frac{1}{2}V(x).
\]
Thus $H(v)$ is sub-$(2^*)$ growth.

The following proposition is important to the computation of a mountain pass level
$c_\varepsilon<\frac{1}{Nm}\varepsilon^nS^{N/2}$.

\begin{prop}\label{prop.H}
Under the assumptions of Theorem \ref{main1},
there exists a function $\tau=\tau(\omega)$ such that
$\lim_{\omega\to0}\tau(\omega)=+\infty$ and for $\omega$ small enough,
\begin{eqnarray*}
\int_{{\mathbb{R}}^N}H(t_\omega v_\omega)\geq\tau(\omega)\cdot\omega^{N-2}.
\end{eqnarray*}
\end{prop}
{\bf Proof}\quad
We divide the proof into three steps.

Step 1: we prove that
\begin{eqnarray}\label{prop.H.1}
\frac{1}{\omega^{N-2}}\int_{B_\omega}H(t_\omega v_\omega)\geq\tau_1(\omega)
\end{eqnarray}
with $\lim_{\omega\to0}\tau_1(\omega)=+\infty$.

By the definition of $v_\omega$, for $x\in B_\omega$, there exist constants $c_2\geq c_1>0$ such that
for $\omega$ small enough, we have
\begin{equation*}\label{prop.H.2}
c_1\omega^{-(N-2)/2}\leq v_\omega(x)\leq c_2\omega^{-(N-2)/2}.
\end{equation*}
and
\begin{equation}\label{prop.H.3}
c_1\omega^{-(N-2)/2}\leq f^m(v_\omega(x))\leq c_2\omega^{-(N-2)/2}.
\end{equation}
On one hand, by (7) of Lemma \ref{lem.f}, (\ref{prop.H.3}) and the continuity of $V(x)$ in $\bar B_\omega$,
there exists $C_1>0$ such that
\begin{equation}\label{prop.H.4}
\int_{B_\omega}V(x)f^2(t_\omega v_\omega)
    \leq C_1\omega^{N-\frac{2}{m}\frac{N-2}{2}}
    = C_1\omega^{(\frac{2^*}{2}-\frac{1}{m})(N-2)}.
\end{equation}
Similarly, there exists $C_2>0$ such that
\begin{equation}\label{prop.H.5}
\int_{B_\omega}f^q(t_\omega v_\omega)
    \geq C_2\omega^{N-\frac{q}{m}\frac{N-2}{2}}
    =C_2\omega^{(\frac{2^*}{2}-\frac{\bar q}{2})(N-2)},
\end{equation}
where $\bar q=q/m$.
On the other hand, using H\"{o}lder inequality, we have
\begin{eqnarray}\label{prop.H.6}
&&\frac{1}{2^*m}\int_{B_\omega}\bigg[(\sqrt{m}t_\omega v_\omega)^{2^*}
        -(f^{m}(t_\omega v_\omega))^{2^*}\bigg] \nonumber\\
    &&\quad\leq\frac{1}{m}\int_{B_\omega}(\sqrt{m}t_\omega v_\omega)^{2^*-1}
        [\sqrt{m}t_\omega v_\omega-f^m(t_\omega v_\omega)] \nonumber\\
    &&\quad \leq\frac{1}{m}\bigg(\int_{B_\omega}(\sqrt{m}t_\omega v_\omega)^{2^*}\bigg)^{(2^*-1)/2^*}
        \bigg(\int_{B_\omega}[\sqrt{m}t_\omega v_\omega-f^m(t_\omega v_\omega)]^{2^*}\bigg)^{1/2^*}.
\end{eqnarray}

{\bf Case 1}: $1< m<2$.
From (\ref{prop.H.6}) and (i) of Lemma \ref{lem.f3}, we obtain that there exists $C_3>0$ such that
\begin{eqnarray}\label{prop.H.7}
&&\frac{1}{2^*m}\int_{B_\omega}\bigg[(\sqrt{m}t_\omega v_\omega)^{2^*}
        -(f^{m}(t_\omega v_\omega))^{2^*}\bigg] \nonumber\\
    &&\quad \leq C_3 \omega^{[N-(\frac{2}{m}-1)\frac{N-2}{2}2^*]\frac{1}{2^*}}
        =C_3\omega^{(1-\frac{1}{m})(N-2)}.
\end{eqnarray}
Combining (\ref{prop.H.4}), (\ref{prop.H.5}) and (\ref{prop.H.7}), we have
\begin{eqnarray*}
&&\frac{1}{\omega^{N-2}}\int_{B_\omega}H(t_\omega v_\omega) \nonumber\\
    &&\quad\geq
-C_1\omega^{(\frac{2^*}{2}-\frac{1}{m}-1)(N-2)}
+C_2\omega^{(\frac{2^*}{2}-\frac{\bar q}{2}-1)(N-2)}
-C_3\omega^{-\frac{1}{m}(N-2)}:=\tau_1(\omega).
\end{eqnarray*}
It is obvious that $\frac{2^*}{2}-\frac{1}{m}-1>-\frac{1}{m}$.
No matter which one of conditions (i)-(iv) in Theorem \ref{main1} holds,
we all have $\frac{2^*}{2}-\frac{\bar q}{2}-1<-\frac{1}{m}$.
It results that $\tau_1(\omega)\to+\infty$ as $\omega\to0$.

{\bf Case 2}: $m\geq2$. Note that for any $\delta\in(0,m)$,
$\lim_{v\to+\infty}\log{f(v)}/f^{\delta}(v)=0,$
we have $\log{f(v)}\leq f^{\delta}(v)$ for $v>0$ large enough.
Thus for $\omega>0$ small enough, from (\ref{prop.H.6}) and (ii) of Lemma \ref{lem.f3}, we get
\begin{eqnarray}\label{prop.H.8}
&&\frac{1}{2^*m}\int_{B_\omega}\bigg[(\sqrt{m}t_\omega v_\omega)^{2^*}
        -(f^{m}(t_\omega v_\omega))^{2^*}\bigg] \nonumber\\
    &&\quad \leq C_3' \omega^{[N-\frac{\delta}{m}\frac{N-2}{2}2^*]\frac{1}{2^*}}
        =C_3'\omega^{\frac{1}{2}(1-\frac{\delta}{m})(N-2)}.
\end{eqnarray}
Combining (\ref{prop.H.4}), (\ref{prop.H.5}) and (\ref{prop.H.8}), we have
\begin{eqnarray*}
&&\frac{1}{\omega^{N-2}}\int_{B_\omega}H(t_\omega v_\omega) \nonumber\\
    &&\quad\geq
-C_1\omega^{(\frac{2^*}{2}-\frac{1}{m}-1)(N-2)}
+C_2\omega^{(\frac{2^*}{2}-\frac{\bar q}{2}-1)(N-2)}
-C_3'\omega^{-\frac{1}{2}(1+\frac{\delta}{m})(N-2)}:=\tau_1(\omega).
\end{eqnarray*}
Since $m\geq2$, we have $\frac{2^*}{2}-\frac{1}{m}-1>-\frac{1}{2}(1+\frac{\delta}{m})$.
No matter which one of conditions (v)-(viii) in Theorem \ref{main1} holds,
there exists a $\delta=\delta(N,\bar q)>0$ (depends on $N$ and $\bar q$) small enough such that
$\frac{2^*}{2}-\frac{\bar q}{2}-1<-\frac{1}{2}(1+\frac{\delta}{m})$.
It results that $\tau_1(\omega)\to+\infty$ as $\omega\to0$.

Case 1 and case 2 show that (\ref{prop.H.1}) holds.

Step 2: we prove that there exists $C_4>0$ such that
\begin{eqnarray}\label{prop.H.9}
\frac{1}{\omega^{N-2}}\int_{B_{2\omega}\setminus B_\omega}H(t_\omega v_\omega)\geq -C_4\omega^{(\frac{2^*}{2}-\frac{1}{m}-1)(N-2)}:=\tau_2(\omega).
\end{eqnarray}

Note that for $x\in B_{2\omega\setminus B_\omega}$, we have
\begin{eqnarray}\label{prop.H.10}
v_\omega(x)\leq v^*_\omega(x) \leq c_2\omega^{-(N-2)/2}.
\end{eqnarray}
Since $\eta_\omega$ is a positive smooth cut-off function,
without lost of generality, we may assume that $\eta_\omega$ is such that
\[
\int_{B_{2\omega}\setminus B_\omega}|v_\omega|^{2^*}
    \leq \int_{B_{2\omega}\setminus B_\omega}V(x)f^2(v_\omega).
\]
Thus by (\ref{lem.f.limit}) and (\ref{prop.H.10}), we have
\begin{eqnarray*}
&&\frac{1}{\omega^{N-2}}\int_{B_{2\omega}\setminus B_\omega}H(t_\omega v_\omega)\\
    &&\quad \geq-\frac{1}{2\omega^{N-2}}\int_{B_{2\omega}\setminus B_\omega}V(x)f^2(t_\omega v_\omega)
        -\frac{1}{2^*m\omega^{N-2}}\int_{B_{2\omega}\setminus B_\omega}|\sqrt{m}t_\omega v_\omega|^{2^*} \\
    &&\quad \geq-\frac{C_5}{\omega^{N-2}}\int_{B_{2\omega}\setminus B_\omega}V(x)f^2(t_\omega v_\omega)\\
    &&\quad \geq-C_4\omega^{N-\frac{2}{m}\frac{N-2}{2}-(N-2)}
        = -C_4\omega^{(\frac{2^*}{2}-\frac{1}{m}-1)(N-2)},
\end{eqnarray*}
where $C_4>0$, $C_5>0$ are constants. This shows that (\ref{prop.H.9}) holds.

Step 3: to conclude, let $\tau(\omega)=\tau_1(\omega)+\tau_2(\omega)$,
we have $\tau(\omega)\to+\infty$ as $\omega\to0$. This implies the conclusion of the proposition.
$\quad\Box$

$\\$
\medskip
{\bf Proof of Theorem \ref{main1}}\quad
By (4) and (7) of Lemma \ref{lem.f}, it is easy to verify that $\bar J$ has the Mountain Pass Geometry.
Lemma \ref{lem.PS} shows that $\bar J$ satisfies (PS) condition.
We prove that $\bar J$ has the mountain pass level $c_\varepsilon<\frac{1}{Nm}\varepsilon^NS^{N/2}$.
Let
\[
F(t)=\frac{\varepsilon^2}{2}\|\nabla(tv_\omega)\|_{L^2}^2-\frac{1}{2^*m}\|\sqrt{m}tv_\omega\|_{L^{2^*}}^{2^*}.
\]
Then we have
\begin{eqnarray}\label{main1.1}
F(t)\leq F(t_0)=\frac{1}{Nm}\varepsilon^NS^{N/2}+O(\omega^{N-2}),\quad\forall t\geq0,
\end{eqnarray}
where $t_0=\frac{\varepsilon^{(N-2)/2}}{\sqrt{m}}$.
By (\ref{main1.1}) and (\ref{prop.H}), we have
\begin{eqnarray*}
\bar J(t_\omega v_\omega)
    &=&F(t_\omega v_\omega)-\int_{{\mathbb{R}}^N}H(t_\omega v_\omega) \\
    &\leq& \frac{1}{Nm}\varepsilon^NS^{N/2}+O(\omega^{N-2})-\tau(\omega)\omega^{N-2} \\
    &<&\frac{1}{Nm}\varepsilon^NS^{N/2}.
\end{eqnarray*}
This shows that $\bar J(v)$ has a nontrivial critical point $v_\varepsilon\in X$, which is a weak solution of (\ref{Sch4}).

We prove that $v_\varepsilon$ is also a weak solution of (\ref{Sch3}).
Firstly, we can argue as the proof of Proposition 2.1 in \cite{Pi-Fe-96} to obtain that
\[
\lim_{\varepsilon\to0}\max_{x\in\partial B_R}{v_\varepsilon(x)}=0.
\]
Thus there exists $\varepsilon_0>0$ such that for all $\varepsilon\in(0,\varepsilon_0)$, we have
$v_\varepsilon(x)\leq a:=f^{-1}(l),~\forall ~|x|= R$, where $l$ is given in (\ref{def.l}).
Secondly, we prove that
\begin{equation}\label{main1.2}
    v_\varepsilon(x)\leq a,
        \quad\forall~\varepsilon\in(0,\varepsilon_0)\mbox{ and }\forall ~x\in{\mathbb{R}}^N\setminus B_R.
\end{equation}
Taking
\[
\varphi=\left\{
          \begin{array}{ll}
            (v_\varepsilon-a)^+, & \hbox{$x\in{\mathbb{R}}^N\setminus B_R$;} \\
            0, & \hbox{$x\in B_R$.}
          \end{array}
        \right.
\]
as a test function in $\langle\bar J'(v_\varepsilon),\varphi\rangle=0$, we get
\begin{eqnarray}\label{main1.3}
&&\varepsilon^2\int_{{\mathbb{R}}^N\setminus B_R}|\nabla (v_\varepsilon-a)^+|^2 \nonumber\\
    &&\quad +\varepsilon^2\int_{{\mathbb{R}}^N\setminus B_R}\bigg(V(x)-\frac{p(x,f(v_\varepsilon))}{f(v_\varepsilon)}\bigg)
f(v_\varepsilon)f'(v_\varepsilon)(v_\varepsilon-a)^+=0.
\end{eqnarray}
By ($p_2$), we have
\[
V(x)-\frac{p(x,f(v_\varepsilon))}{f(v_\varepsilon)}>0,\quad\forall x\in{\mathbb{R}}^N\setminus B_R.
\]
Therefore, all terms in (\ref{main1.3}) must be equal to zero. This implies $v_\varepsilon\leq a$
in ${\mathbb{R}}^N\setminus B_R$.
This proves (\ref{main1.2}). Thus $v_\varepsilon$ is a solution of Problem (\ref{Sch3}).

To complete the proof, we deduce as the proof for Theorem 4.1 in \cite{Li-11}
to obtain that $v_\varepsilon|_{B_R}\in L^{\infty}(B_R)$. Thus $u_\varepsilon=f(v_\varepsilon)\in X\cap L^{\infty}({\mathbb{R}}^N)$ is a nontrivial weak solution of (\ref{sch2}).
Finally, by Proposition \ref{prop.u} in the following, we have
$
\lim_{\varepsilon\to0}\|u_\varepsilon\|_X=0
$
and
$
u_\varepsilon(x)\leq Ce^{-\frac{\beta}{\varepsilon}|x-x_\varepsilon|}.
$
This completes the proof.
$\quad\Box$

We prove the norm estimate and the exponential decay.

\begin{prop}\label{prop.u}
Let $v_\varepsilon\in X\cap L^{\infty}({\mathbb{R}}^N)$ be a solution of (\ref{Sch3})
and let $u_\varepsilon=f(v_\varepsilon)$, then we have
\[
\lim_{\varepsilon\to0}\|u_\varepsilon\|_X=0,\quad\mbox{and}\quad
u_\varepsilon(x)\leq Ce^{-\frac{\beta}{\varepsilon}|x-x_\varepsilon|},
\]
where $C>0$, $\beta>0$ are constants.
\end{prop}
{\bf Proof}\quad
Firstly, let $x_0\in B_R$ be such that $V(x_0)=V_0$. Define $J_{0}: X\to {\mathbb{R}}$ by
\[
J_{0}(v)=\frac{1}{2}\int_{{\mathbb{R}}^N}|\nabla v|^2+\frac{1}{2}\int_{{\mathbb{R}}^N}V_0f^2(v)
-\int_{{\mathbb{R}}^N}G(f(v)).
\]
Let
\[
c_0=\inf_{\gamma\in\Gamma_0}\sup_{t\in[0,1]}J_{0}(\gamma(t)),
\]
\[
\Gamma_0=\{v\in C([0,1],X):\gamma(0)=0, J_{0}(\gamma(1))<0\}.
\]
Similar to the proof for estimate (2.4) in \cite{Pi-Fe-96} (or Lemma 3.1 in \cite{Sq-03}),
we can show that
$
c_\varepsilon\leq \varepsilon^Nc_0+o(\varepsilon^N)
$
by using the change of coordinates $y=(x-x_0)/\varepsilon$.
Arguing as for (\ref{lem2-6-0}), and by virtue of this energy estimate, we obtain
\[
\|v_\varepsilon\|_X\leq\frac{\theta c_\varepsilon}{\min\{(\frac{\theta}{2}-m)\varepsilon^2,(\frac{\theta}{2}-\frac{\theta}{2k}-1)\}}
\leq\frac{2\theta c_0}{\theta-2m}\varepsilon^{N-2}+o(\varepsilon^{N-2})
\]
for $\varepsilon>0$ sufficient small.
Let $u_\varepsilon=f(v_\varepsilon)$, then $u_\varepsilon\neq0$.
Note that $|\nabla u_\varepsilon|\leq |v_\varepsilon|$ and $|u_\varepsilon|\leq|v_\varepsilon|$, we get
$
\lim_{\varepsilon\to0}\|u_\varepsilon\|_X=0.
$

Secondly, similar to the proof for Theorem 4.1 in \cite{Li-11}, we conclude that $v_\varepsilon\in L^{\infty}({\mathbb{R}}^N)$ and by \cite{La-Ur-68}, we have $v_\varepsilon\in C^{1,\alpha}(B_R)$.
Now let $x_\varepsilon$ denote the maximum point of $v_\varepsilon$ in $B_R$ and let
\[
\sigma:=\sup\{s>0:g(t)< V_0t \mbox{ for every } t\in[0,s]\}.
\]
Then $v_\varepsilon(x_\varepsilon)\geq f^{-1}(\sigma)$ for $\varepsilon>0$ small.
In fact, assume that $v_\varepsilon(x_\varepsilon)<f^{-1}(\sigma)$ for some $\varepsilon>0$ sufficiently small.
According to the definition of $l$ (see (\ref{def.l})) and $\sigma$, we have
$v_\varepsilon(x)\leq f^{-1}(l)<f^{-1}(\sigma)$ (note that $k>1$ in (\ref{def.l})),
$\forall x\in {\mathbb{R}}^N\setminus B_R$. Thus
\[
V(x)-\frac{g(f(v_\varepsilon))}{f(v_\varepsilon)}>0,\quad\forall x\in{\mathbb{R}}^N.
\]
Since $v_\varepsilon=f^{-1}(u_\varepsilon)$ is a critical point of $J_\varepsilon$,
we choose $\varphi=f(v_\varepsilon)/f'(v_\varepsilon)$
as a test function in $\langle J'_\varepsilon(v_\varepsilon),\varphi\rangle=0$ and get
\begin{eqnarray*}
0&=& \varepsilon^2\int_{{\mathbb{R}}^N}\bigg(1+\frac{m(m-1)|f(v_n)|^{2(m-1)}}{1+m|f(v_n)|^{2(m-1)}}\bigg)|\nabla v_\varepsilon|^2 \\
        &&+\int_{{\mathbb{R}}^N}V(x)f^2(v_\varepsilon)-\int_{{\mathbb{R}}^N}g(f(v_\varepsilon))f(v_\varepsilon) \\
&=& \varepsilon^2\int_{{\mathbb{R}}^N}\bigg(1+\frac{m(m-1)|f(v_n)|^{2(m-1)}}{1+m|f(v_n)|^{2(m-1)}}\bigg)|\nabla v_\varepsilon|^2 \\
        &&+\int_{{\mathbb{R}}^N}\bigg(V(x)-\frac{g(f(v_\varepsilon))}{f(v_\varepsilon)}\bigg)f^2(v_\varepsilon).
\end{eqnarray*}
It turns out that all terms in the above equality must be equal to zero, which means that $v_\varepsilon\equiv0$, a contradiction.

Now let $w_\varepsilon(x)=v_{\varepsilon}(x_\varepsilon+\varepsilon x)$, then $w_\varepsilon$ solves the equation
\begin{eqnarray*}
    -\Delta w_\varepsilon+V(x_\varepsilon+\varepsilon x)f(w_\varepsilon)f'(w_\varepsilon)=g(f(w_\varepsilon))f'(w_\varepsilon), \quad x\in{\mathbb{R}}^N.
\end{eqnarray*}
Note that $\lim_{t\to0^+}\frac{f(t)f'(t)}{t}=1$ by the properties of $f$ and that
$w_\varepsilon(x)\to0$ as $|x|\to+\infty$, we have,
there exists $R_0>0$ such that for all $|x|\geq R_0$,
\begin{equation}\label{thm-3}
    f(w_\varepsilon(x))f'(w_\varepsilon(x))\geq\frac{3}{4}w_\varepsilon(x)
\end{equation}
and
\begin{equation}\label{thm-4}
    g(f(w_\varepsilon(x)))f'(w_\varepsilon(x))\leq\frac{V_0}{2}w_\varepsilon(x).
\end{equation}
Let $\varphi(x)=Me^{-\beta|x|}$ with $\beta^2<\frac{V_0}{4}$ and
$Me^{-\beta R_0}\geq w_\varepsilon(x)$ for all $|x|=R_0$. It is easy to verify that for $x\neq0$,
\begin{equation}\label{thm-5}
    \Delta \varphi\leq\beta^2\varphi.
\end{equation}
Now define $\psi_\varepsilon=\varphi-w_\varepsilon$. Using (\ref{thm-3})-(\ref{thm-5}), we have
\begin{eqnarray*}
\left\{
  \begin{array}{ll}
    -\Delta \psi_\varepsilon+\frac{V_0}{4}\psi_\varepsilon\geq0, & \hbox{in $|x|\geq R_0$;} \\
    \psi_\varepsilon\geq0, & \hbox{in $|x|=R_0$;} \\
    \lim_{|x|\to\infty}\psi_\varepsilon=0. & \hbox{}
  \end{array}
\right.
\end{eqnarray*}
By the maximum principle, we have $\psi_\varepsilon\geq0$ for all $|x|\geq R_0$. Thus, we obtain that for all $|x|\geq R_0$,
\[
w_\varepsilon(x)\leq\varphi(x)\leq Me^{-\beta|x|}.
\]
Using the change of variable, we have that for all $|x|\geq R_0$,
\[
v_\varepsilon(x)=w_\varepsilon(\varepsilon^{-1}(x-x_\varepsilon))\leq Me^{-\frac{\beta}{\varepsilon}|x-x_\varepsilon|}.
\]
Then by the regularity of $v_\varepsilon$ on $B_R$ and note that $f(t)\leq t$ for all $t\geq0$, we have
\[
u_\varepsilon(x)\leq Ce^{-\frac{\beta}{\varepsilon}|x-x_\varepsilon|}
\]
for some $C>0$. This completes the proof.
$\quad\square$

$\\$
\medskip
{\bf Proof of Theorem \ref{main2}}\quad
We consider the following equation
\begin{eqnarray}\label{sch5}
-\Delta u&+&V(x)u-k\alpha(\Delta(|u|^{2\alpha}))|u|^{2\alpha-2}u \nonumber\\
    &=&\lambda|u|^{q-2}u+|u|^{2^*(2\alpha)-2}u,\quad u>0,x\in{\mathbb{R}}^N.
\end{eqnarray}
Let $y=\varepsilon x$ with $\varepsilon\in(0,\varepsilon_0)$, $\varepsilon_0$ is given by Theorem \ref{main1},
then we can transform (\ref{sch5}) into
\begin{eqnarray}\label{sch6}
-\varepsilon^2\Delta u&+&\bar V(y)u-k\alpha\varepsilon^2(\Delta(|u|^{2\alpha}))|u|^{2\alpha-2}u \nonumber\\
    &=&\lambda|u|^{q-2}u+|u|^{2^*(2\alpha)-2}u,\quad u>0,y\in{\mathbb{R}}^N.
\end{eqnarray}
Here $\bar V(y)=V(\frac{y}{\varepsilon})$ still has the properties given in assumption (V).
Thus according to Theorem \ref{main1}, (\ref{sch6}) has a positive weak solution
$u_\varepsilon(y)$ in $X\cap L^\infty({\mathbb{R}}^N)$, this implies that
(\ref{sch5}) has a positive weak solution $u_1(x)=u_\varepsilon(\varepsilon x)$.
$\quad\Box$

%\begin{acknowledgements}
%\end{acknowledgements}



\end{document}